\documentclass[a4paper]{amsart}

\usepackage[T1]{fontenc}
\usepackage[utf8]{inputenc}
\usepackage{textcomp}
\usepackage{lmodern}
\usepackage{amsmath,amsthm,amssymb}
\usepackage{microtype}
\usepackage{graphicx}
\usepackage{todonotes}
\usepackage{commath}
\usepackage[numbers,sort,compress]{natbib}
\usepackage{hyperref}
\hypersetup{%
  pdfauthor={Matthias. J. Weber, Hans-Peter Schröcker},%
  pdftitle={Minimal area ellipses in the hyperbolic plane},%
  pdfkeywords={hyperbolic geometry, enclosing ellipse, minimal area, uniqueness},%
  colorlinks=true,%
  citecolor=black,%
  urlcolor=black,%
  linkcolor=black,%
}

\newcommand*{\HP}{\mathbb{H}^2}
\newcommand*{\US}{\mathbb{S}^2_1}
\newcommand*{\RSet}{\mathbb{R}}
\newcommand*{\mNorm}[1]{\Vert #1 \Vert}
\newcommand*{\mIP}[2]{\langle #1, #2 \rangle}
\newcommand*{\tp}{\mathrm{T}}
\newcommand*{\ellE}{\bar{E}}
\newcommand*{\ellK}{\bar{K}}
\newcommand*{\cA}{\mathrm{A}}
\newcommand*{\cB}{\mathrm{B}}
\newcommand*{\cC}{\Gamma}
\newcommand*{\cD}{\Delta}
\newcommand*{\conv}[1]{\overline{#1}}
\DeclareMathOperator{\dist}{dist}
\DeclareMathOperator{\diag}{diag}
\DeclareMathOperator{\area}{area}
\DeclareMathOperator{\arccosh}{arccosh}

\newtheorem{theorem}{Theorem}
\newtheorem{lemma}[theorem]{Lemma}
\newtheorem{proposition}[theorem]{Proposition}

\theoremstyle{definition}
\newtheorem{definition}[theorem]{Definition}
\newtheorem{example}[theorem]{Example}

\begin{document}

\title{Minimal area ellipses in the hyperbolic plane}
\author{Matthias.\,J.~Weber \and Hans-Peter~Schröcker}
\date{\today}

\address{Matthias.\,J.~Weber, Hans-Peter Schröcker, Unit Geometry and
  CAD, University Innsbruck, Technikerstraße 13, 6020 Innsbruck,
  Austria}

 % hyperbolic and elliptic geometries (general) and
         % generalizations

\begin{abstract}
  We present uniqueness results for enclosing ellipses of minimal area
  in the hyperbolic plane. Uniqueness can be guaranteed if the
  minimizers are sought among all ellipses with prescribed axes or
  center. In the general case, we present a sufficient and easily
  verifiable criterion on the enclosed set that ensures uniqueness.
\end{abstract}

\maketitle

\section{Introduction and statement of the main result}
\label{sec:introduction}

By a well-known theorem of convex geometry, a full-dimensional,
compact subset $F$ of the Euclidean plane can be enclosed by a unique
ellipse $C$ of minimal area.  We share the general belief that this is
an important but easy result.  Therefore, it is not surprising that
recently much more general uniqueness results were obtained
\cite{gruber08:_john_type,%
  schroecker08:_uniqueness_results_ellipsoids,%
  weber10:_davis_convexity_theorem}. These articles also contain more
complete references to the relevant literature.

The situation in the elliptic plane is different: As of today,
uniqueness can only be guaranteed for ``sufficiently small and round
sets $F$''. The precise statement can be found in
\cite[Theorem~8]{weber10:_minimal_area_conics}.  Its proof requires
some non-trivial calculations.  It is still an open question whether
any compact subset of the elliptic plane possesses a unique enclosing
conic or not.

In this article we consider uniqueness of the minimal area ellipse in
the hyperbolic plane.  The algebraic equivalence of elliptic and
hyperbolic geometries suggests to imitate the proof of
\cite[Theorem~8]{weber10:_minimal_area_conics}.  Indeed, this is
possible to a large extent, but not completely.  The outcome of this
research is similar to the elliptic case.  Uniqueness can be
guaranteed if some conditions on the axis lengths of enclosing
ellipses of minimal area are met.  If this is not possible, we can
make neither a positive nor a negative uniqueness statement.  Our main
result is

\begin{theorem}
  \label{th:1}
  Consider a compact and full-dimensional subset $F$ of the hyperbolic
  plane. The enclosing ellipse of minimal area to $F$ is unique if the
  following conditions are met:
  \begin{itemize}
  \item There exist positive numbers $\varrho$, $\mathrm{R}$ such that
    the semi-axis lengths of the (a priori not necessarily unique)
    minimal ellipses are in the closed interval
    $[\varrho,\mathrm{R}]$.
  \item The values $\nu_1 = \coth^2\!\mathrm{R}$ and $\nu_2 =
    \coth^2\!\varrho$ satisfy the inequality
    \begin{equation}
      \label{eq:1}
      H(\nu_1,\nu_2) := -13\nu_1^2 + 5\nu_1\nu_2 - 3\nu_1 + 7\nu_2 + 4 \le 0.
    \end{equation}
  \end{itemize}
\end{theorem}

The minimal enclosing ellipse $C_{\min}$ to the convex hull $F$ of a
finite point set is depicted in Figure~\ref{fig:center-and-axes}.  The
drawing refers to the Cayley-Klein model of the hyperbolic plane which
will be introduced in Section~\ref{sec:hyperboloid-model}.

\begin{figure}
  \centering
  \includegraphics{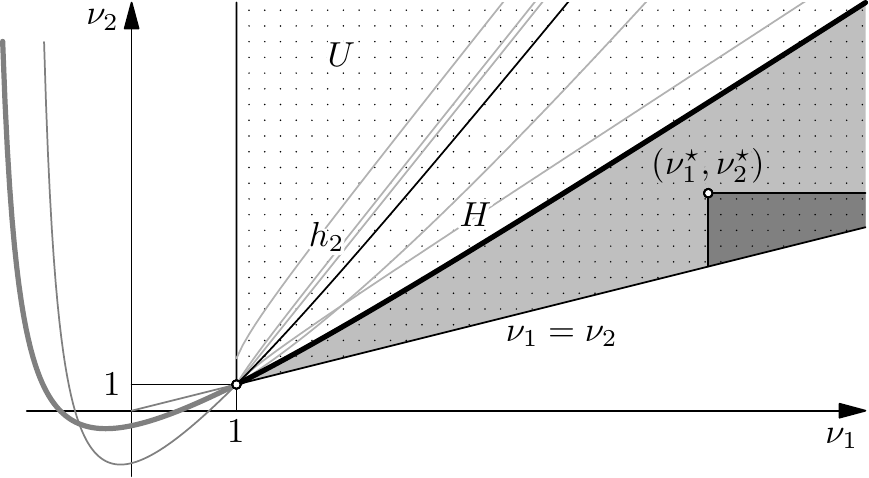}
  \caption{The curves $H(\nu_1,\nu_2) = 0$ and $h_2(\nu_1,\nu_2) = 0$.}
  \label{fig:hyp-graph}
\end{figure}

Figure~\ref{fig:hyp-graph} depicts the graph of the function
$H(\nu_1,\nu_2) = 0$.  The shaded area contains admissible values
$\nu_1$, $\nu_2$.  The meaning of the remaining elements will be
explained later in the text.

Condition \eqref{eq:1} is certainly fulfilled for $\nu_1 = \nu_2$ and
$\nu_1 \to \infty$. Thus, Theorem~\ref{th:1} informally states that
minimal enclosing ellipses are unique if they are sufficiently small
and round.  Of course, the Theorem should be accompanied by an easily
verifiable criterion that ensure suitably shaped minimal ellipses:

\begin{proposition}
  \label{prop:1}
  Consider a compact and full-dimensional subset $F$ of the hyperbolic
  plane and denote its (hyperbolic) convex hull by $\conv{F}$. Assume
  $\conv{F}$ admits an inscribed circle of radius $\varrho$ and a
  circumscribed ellipse of area $S$. Denote by $\mathrm{R}$ the major
  semi-axis length of an ellipse of area $S$ and minor semi-axis
  length $r$. Then the minimal area ellipse of $F$ has semi-axis
  length in the interval $[\varrho,\mathrm{R}]$.
\end{proposition}

We omit the obvious proof of this proposition.  Together with
Theorem~\ref{th:1}, it leads to the following sufficient test for the
uniqueness of the minimal enclosing ellipsoid to a given set~$F$:
\begin{enumerate}
\item Find a (large) inscribed circle to $\conv{F}$ and denote its
  radius by~$\varrho$.
\item Find a (small) circumscribed ellipse to $\conv{F}$ and denote
  its area by~$S$.
\item Compute the unique value $\mathrm{R}$ such that an ellipse with
  semi-axis lengths $\varrho$ and $\mathrm{R}$ has area $S$. By
  construction, $[\varrho,\mathrm{R}]$ is not empty.
\item The minimal area ellipse is unique, if $\varrho$ and
  $\mathrm{R}$ satisfy the inequality \eqref{eq:1}.
\end{enumerate}

We will show later that if a point $(\nu_1,\nu_2)$ satisfies
\eqref{eq:1} than the same is true for every admissible point
$(\nu'_1,\nu'_2)$ with $\nu'_1 \ge \nu_1$ and $\nu'_2 \le \nu_2$. This
means, that chances for an affirmative uniqueness statement increase
with a large value of the radius~$\varrho$ and a small value of the
area $S$, that is, with the quality of the input obtained from the
first and the second step.

The basic ideas and the initial calculations in our proof of
Theorem~\ref{th:1} are more or less identical to the proof of
\cite[Theorem~8]{weber10:_minimal_area_conics}.  The minor differences
pertain to occasional changes in sign and the use of the hyperbolic
functions $\cosh$, $\sinh$, etc. instead of their spherical
counterparts $\cos$, $\sin$, etc.  The major differences are in the
final estimates.  Given the similarities between the elliptic and
hyperbolic case, we consider a rather terse presentation
appropriate. Yet, we will try to work out the crucial junctions points
and the major differences.

In Section~\ref{sec:hyperboloid-model} we settle our notation and
introduce the hyperboloid model of the hyperbolic plane, where our
calculations take place. In Section~\ref{sec:area-ellipses} we provide
a formula for the area of ellipses in the hyperbolic plane, which is
probably hard to find elsewhere.  In Section~\ref{sec:area-convex}, we
prove a fundamental convexity result for the area function.  By
standard arguments, it yields uniqueness of the minimal ellipse among
all ellipses with prescribed axes or center. The proof of
Theorem~\ref{th:1} is given in Section~\ref{sec:uniqueness}.  Its main
ingredient is Lemma~\ref{lem:half-turn}, the Half-Turn Lemma. The
merely technical parts of its proof are moved to the appendix.

\section{The hyperboloid model of hyperbolic geometry}
\label{sec:hyperboloid-model}

In \cite{weber10:_minimal_area_conics} we used the spherical model of
the elliptic plane for investigating uniqueness of minimal area
conics. It is obtained from the geometry of the unit sphere
$\mathbb{S}^2$ of Euclidean three-space by identifying antipodal
points. By analogy, our calculations in this article refer to the
spherical model (or ``hyperboloid model'') of the hyperbolic plane
which is obtained in similar fashion from the geometry of the sphere
of squared radius $-1$ in Minkowski three space $\RSet^3_1$. An
elementary introduction to this model is given in
\cite{reynold93:_hyperbolic_geometry}.

Minkowski three-space $\RSet^3_1$ is the metric space over $\RSet^3$
where the metric is induced by the indefinite inner product
\begin{equation*}
  \mIP{x}{y} = -x_0y_0 + x_1y_1 + x_2y_2.
\end{equation*}
The locus of the spherical model of the hyperbolic plane is the sphere
$\US$, defined as
\begin{equation*}
  \US = \{x \in \RSet^3_1\colon \mNorm{x}^2 = -x_0^2 + x_1^2 + x_2^2 = -1\}.
\end{equation*}
In a Euclidean interpretation, it is a hyperboloid of two sheets. We
use $\US$ as a model of the hyperbolic plane $\HP$. The following
concepts are taken from \cite{reynold93:_hyperbolic_geometry}:
\begin{itemize}
\item The points of $\HP$ are the points of $\US$ with antipodal
  points $x$ and $-x$ identified.
\item The lines of $\HP$ are the intersections of $\US$ with planes
  through the origin $0$.
\item The hyperbolic distance between two points $x$, $y \in \US$ is
  defined by
  \begin{equation*}
    \dist(x,y) = \arccosh(-\mIP{x}{y}).
  \end{equation*}
\item The hyperbolic angle between two straight lines $K$ and $L$ is
  defined by
  \begin{equation*}
    \sphericalangle(k,l) = \arccos\frac{\mIP{k}{l}}{\mNorm{k} \cdot
      \mNorm{l}}
  \end{equation*}
  where $k$ and $l$ are two arbitrary tangent vectors of $K$ and $L$,
  respectively.
\end{itemize}

Note that this model of $\HP$ is closely related to the well-known
bundle model and also the Cayley-Klein model of the hyperbolic
plane. The bundle model is obtained by connecting points and lines
from the spherical model with the origin $0$ of $\RSet^3_1$; the
Cayley-Klein model is obtained by intersecting the bundle model with
the plane $x_0 = 1$.  Its points are the inner points of the circle
\begin{equation*}
  K\colon x_0 = 1,\ x_1^2 + x_2^2 = 1.
\end{equation*}
We will occasionally use the Cayley-Klein model for the purpose of
visualization but it is also convenient for defining center and axes
of an ellipse $C$ in the hyperbolic plane.

\begin{figure}
  \centering
  \includegraphics{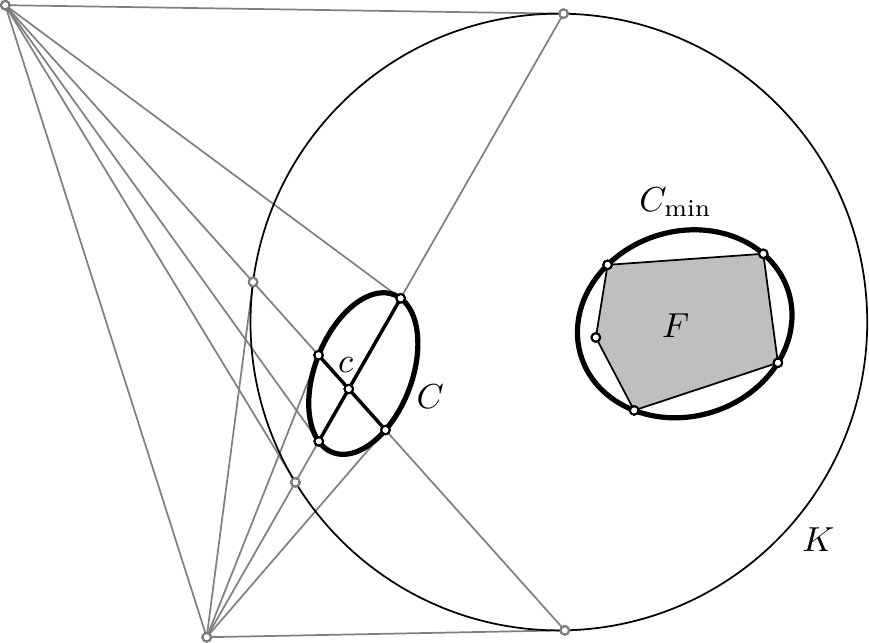}
  \caption{Center $c$ and axes of a hyperbolic ellipse $C$; minimal
    ellipse $C_{\min}$ to the convex hull $F$ of a finite points set}
  \label{fig:center-and-axes}
\end{figure}

The conics in the spherical model of $\HP$ are the intersections of
$\US$ with quadratic cones centered at $0$. In the Cayley-Klein model,
hyperbolic ellipses are conics that lie in the interior of $K$.  The
ellipse center is the unique vertex $c$ of the common polar triangle
$P$ of $C$ and $K$. It is indeed a center in elementary sense, as it
halves the (hyperbolic) distance between the ellipse points on any line
incident with~$c$. The axes of $C$ are the two sides of $P$ through
$c$. Degenerate pole triangles characterize the circles among the
ellipses. Their center is still well-defined but the axes are
undetermined so that any line through $c$ can be addressed as axis.
Figure~\ref{fig:center-and-axes} displays a hyperbolic ellipse, its
center and axes in the Cayley-Klein model.

\section{The area of ellipses}
\label{sec:area-ellipses}

The hyperbolic plane $\HP$ can be parametrized as
\begin{equation}
  \label{eq:2}
  \HP\colon Y(\theta,\varphi) =
  \begin{pmatrix}
    \cosh\theta\\
    \sinh\theta \sin\varphi\\
    \sinh\theta \cos\varphi
  \end{pmatrix},
  \quad \theta \in [0, \infty), \varphi \in [-\pi, \pi).
\end{equation}
A conic $C$, defined as the intersection of this point set with a
quadratic cone whose vertex is in the origin, can be described as
\begin{equation*}
  C = \{x \in \HP \colon x^\tp \cdot M \cdot x = 0\},
\end{equation*}
where $M \in \RSet^{3 \times 3}$ is an indefinite symmetric matrix of
full rank. A vector $x$ is called (Minkowski) eigenvector of $M$ with
(Minkowski) eigenvalue $\lambda$ if
\begin{equation}
  \label{eq:3}
  M \cdot x = \lambda I \cdot x.
  \quad\text{where}\quad
  I = \diag(-1, 1, 1).
\end{equation}
By $e(M) = (\nu_0,\nu_1,\nu_2)$ we denote the vector of eigenvalues of
$M$, arranged in ascending order. We will only consider the case where
$M$ describes an ellipse. In this case $M$ can be normalized such that
$e(M) = (1, \nu_1, \nu_2)$ and $1 < \nu_1 \le \nu_2$.

A point $x$ is contained in the ellipse $C$ if it satisfies $x^\tp
\cdot M \cdot x < 0$ and $M$ is in normal form. After a suitable
(Minkowski) rotation of $\US$ we may assume that the ellipse is
described by the diagonal matrix
\begin{equation}
  \label{eq:4}
  M = \diag(-1, \nu_1, \nu_2).
\end{equation}
Referring to the parametrization \eqref{eq:2}, points inside $C$
belong to parameter values $(\theta, \varphi)$ related by
\begin{equation*}
  \theta  < \theta^\star =
  \arccosh \sqrt{\frac{\nu_1 \sin^2\!\varphi + \nu_2 \cos^2\!\varphi}
                      {\nu_1 \sin^2\!\varphi + \nu_2 \cos^2\!\varphi - 1}}.
\end{equation*}
By integrating the area element
\begin{equation*}
  \sqrt{
    \Bigl\langle \dpd{H}{\theta}, \dpd{H}{\theta} \Bigr\rangle
    \cdot
    \Bigl\langle \dpd{H}{\varphi}, \dpd{H}{\varphi} \Bigr\rangle
    -
    \Bigl\langle \dpd{H}{\theta}, \dpd{H}{\varphi} \Bigr\rangle^2
  } \dif\theta \wedge \dif\varphi =
  \sinh\theta \dif\theta \wedge \dif\varphi
\end{equation*}
of \eqref{eq:2} (see for example Proposition~5.2 of
\cite{callahan00:_geometry_spacetime}) we obtain the area of the conic
$C$ as
\begin{equation}
  \label{eq:5}
  \begin{aligned}
  \area(C) &= \area(\nu_1, \nu_2)
            = \int_{-\pi}^\pi \int_0^{\theta^\star} \sinh\theta \dif\theta \dif\varphi \\
           &= \int_{-\pi}^\pi (\cosh\theta^\star - 1) \dif\varphi
            = \int_{-\pi}^\pi \sqrt{\frac{\nu_1\sin^2\varphi + \nu_2\cos^2\varphi}
                                       {\nu_1\sin^2\varphi + \nu_2\cos^2\varphi - 1}} \dif \varphi - 2 \pi.
  \end{aligned}
\end{equation}
This is valid as long as $M$ is normalized such that $e(M) =
(1,\nu_1,\nu_2)$. If $M$ is not normalized and has ordered eigenvalues
$e(M) = (\nu_0, \nu_1, \nu_2)$, the area formula becomes
\begin{equation}
  \label{eq:6}
  \area(\nu_0,\nu_1,\nu_2) =
  \int_{-\pi}^\pi \sqrt{
        \frac{\nu_1\sin^2\!\varphi + \nu_2\cos^2\!\varphi}
             {\nu_1\sin^2\!\varphi + \nu_2\cos^2\!\varphi - \nu_0}
  }\dif\varphi - 2\pi.
\end{equation}

\section{Convexity of the area function.}
\label{sec:area-convex}

Convexity of the area function \eqref{eq:5} is already the key
property for uniqueness of the minimal area ellipse among concentric
or co-axial ellipses.  Recall that only values $\nu_1$, $\nu_2 > 1$
are admissible.

\begin{lemma}
  \label{lem:1}
  The area function~\eqref{eq:5} is strictly convex for $\nu_1$,
  $\nu_2 > 1$.
\end{lemma}

\begin{proof}
  We proof that the Hessian matrix of~\eqref{eq:5} is positive
  definite, that is, all its principal minors are positive.  The upper
  left entry equals
  \begin{equation}
    \label{eq:7}
    \dpd[2]{\area}{\nu_1} = \frac{1}{4} \int_{-\pi}^\pi J \sin^4\!\varphi \dif \varphi,
  \end{equation}
  where
  \begin{equation*}
    J = \frac{4\nu_1\sin^2\!\varphi + 4\nu_2\cos^2\!\varphi - 1}
             {(\nu_1\sin^2\!\varphi +  \nu_2\cos^2\!\varphi)^{3/2}
              (\nu_1\sin^2\!\varphi +  \nu_2\cos^2\!\varphi - 1)^{5/2}}.
  \end{equation*}
  Clearly, $J$ is positive for admissible values of $\nu_1$ and
  $\nu_2$. Therefore, \eqref{eq:7} is positive as well. The
  determinant of the Hessian matrix is
  \begin{multline}
    \label{eq:8}
    \dpd[2]{\area}{\nu_1} \dpd[2]{\area}{\nu_2} - \Bigl(
      \dmd{\area}{}{\nu_1}{}{\nu_2}{} \Bigr)^2 \\
    =\frac{1}{16} \int_{-\pi}^\pi J \sin^4\!\varphi \dif \varphi \cdot
                 \int_{-\pi}^\pi J \cos^4\!\varphi \dif \varphi -
    \frac{1}{16}
    \Bigl(
        \int_{-\pi}^\pi J \sin^2\!\varphi \cos^2\!\varphi \dif \varphi
    \Bigr)^2.
  \end{multline}
  Because $\sqrt{J}\sin^2\!\varphi$ and $\sqrt{J}\cos^2\!\varphi$ are
  not proportional we can apply the strict Schwarz inequality and find
  \begin{equation*}
    \sqrt{\int_{-\pi}^\pi (\sqrt{J} \sin^2\!\varphi)^2 \dif \varphi}
    \cdot
    \sqrt{\int_{-\pi}^\pi (\sqrt{J} \cos^2\!\varphi)^2 \dif \varphi}
    >
    \int_{-\pi}^\pi J \sin^2\!\varphi \cos^2\!\varphi \dif \varphi.
  \end{equation*}
  Thus, \eqref{eq:8} is positive and $\area(\nu_1,\nu_2)$ is indeed a
  strictly convex function.
\end{proof}

Now, two uniqueness results follow from standard arguments (see
\cite{schroecker08:_uniqueness_results_ellipsoids,%
  weber10:_davis_convexity_theorem} and in particular
\cite{weber10:_minimal_area_conics}).

\begin{theorem}
  \label{th:2}
  Let $F$ be a compact and full-dimensional subset of the hyperbolic
  plane. Among all ellipses with two given axes that contain $F$ there
  exists exactly one with minimal area.
\end{theorem}

\begin{theorem}
  \label{th:3}
  Let $F$ be a compact and full-dimensional subset of the hyperbolic
  plane. Among all ellipses with given center that contain $F$ there
  exists exactly one with minimal area.
\end{theorem}

We give a quick outline of the proofs of Theorem~\ref{th:2} and
\ref{th:3}, mainly because this gives us the opportunity to introduce
an important concept that will be required later.

\begin{definition}[in-between ellipse]
  \label{def:in-between-ellipse}
  Let $C_0$ and $C_1$ be two ellipses
  \begin{equation*}
    C_i = \{x \in \HP \colon x^\tp \cdot M_i \cdot x = 0\}, \quad i=0,1
  \end{equation*}
  where the matrices $M_i$ are indefinite and have Minkowski
  eigenvalues $\nu_{i,0} = 1$ and $\nu_{i,1}$, $\nu_{i,2} > 1$. For
  $\lambda \in (0,1)$, the \emph{in-between ellipse} $C_\lambda$ of
  $C_0$ and $C_1$ is defined as
  \begin{equation*}
    C_\lambda = \{x \in \HP \colon x^\tp \cdot M_\lambda \cdot x\},
  \end{equation*}
  where
  \begin{equation*}
    M_\lambda = (1-\lambda) M_0 + \lambda M_1.
  \end{equation*}
  We also write $C_\lambda = (1-\lambda) C_0 + \lambda C_1$.
\end{definition}

It is obvious that $C_\lambda$ contains the common interior of $C_0$
and $C_1$ and is an ellipse if this interior is not empty. Moreover,
it follows from Lemma~\ref{lem:1} and the strict version of Davis'
convexity theorem \cite{davis57:_convex_functions,%
  lewis96:_convex_analysis} that $\area(C_\lambda)$ is a strictly
convex function of $\lambda$.  More detailed arguments can be found in
\cite{weber10:_davis_convexity_theorem,%
  weber10:_minimal_area_conics}.  The important fact to remember is
that two enclosing conics $C_0$ and $C_1$ of the same area give rise
to an in-between conic $C_\lambda$ of lesser area. Thus, the
assumption of two minimal area conics leads to a contradiction.  Note
that convexity of $\area(C_\lambda)$ in the general (non-concentric)
case is not implied by Davis' convexity theorem.

\section{Uniqueness in the general case.}
\label{sec:uniqueness}

Now we come to the proof of Theorem~\ref{th:1}, the general uniqueness
result.  As usual, existence follows from compactness arguments.  The
basic ideas and initial steps in the proof of uniqueness are not
different from the proof of Theorem~8 in
\cite{weber10:_minimal_area_conics}.  We give an outline:
\begin{itemize}
\item Assume existence of two minimal enclosing ellipses $C_0$
  and~$C_1$.
\item Find the unique (hyperbolic) half-turn $\eta$ (an idempotent
  hyperbolic rotation) such that $C^\star_1 = \eta(C_1)$ and
  $C^\star_1$ is concentric with~$C_0$.
\item Define in-between ellipses $C_\lambda = (1-\lambda) C_0 +
  \lambda C_1$ and $C^\star_\lambda = (1-\lambda) C_0 + \lambda
  C^\star_1$ according to Definition~\ref{def:in-between-ellipse}.
\item Show that there exists $\varepsilon > 0$ such that
  $\area(C_\lambda) < \area(C^\star_\lambda)$ for $0 < \lambda <
  \varepsilon$. Because of $\area(C^\star_\lambda) \le \area(C_0) =
  \area(C_1)$ (with equality iff $C_1^\star = C_0$) this contradicts
  the assumed minimality of $C_0$~and~$C_1$.
\end{itemize}
Existence of $\varepsilon$ in the last step of this program can be
proved by showing the inequality
\begin{equation}
  \label{eq:9}
  \dpd{\area(C^\star_\lambda)}{\lambda}\Big|_{\lambda=0}
  <
  \dpd{\area(C_\lambda)}{\lambda}\Big|_{\lambda=0}.
\end{equation}
The advantage of this approach is that both sides of \eqref{eq:9} can
be readily computed from the normalized equations that describe $C_0$,
$C_1$, and $C^\star_1$. In particular, the cubic problem of
calculating the eigenvalues of the matrices describing $C_\lambda$ or
$C^\star_\lambda$ is avoided.

\begin{figure}
  \centering
  \includegraphics{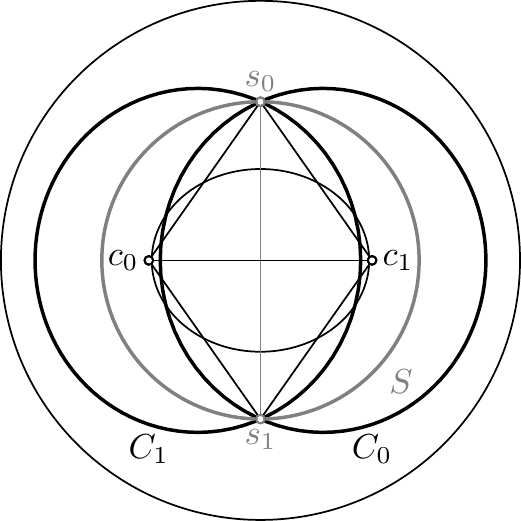}
  \caption{The case of two circles}
  \label{fig:two-circles}
\end{figure}

In order to follow the outline of the proof of Theorem~\ref{th:1} we
have to compute the ellipses $C_0$, $C_1$ and $C^\star_1$ in a
sufficiently general way. By Theorem~\ref{th:3}, the centers of $C_0$
and $C_1$ can be assumed to be different. Thus, there exists a unique
mid-point $r$ of their respective centers $c_0$ and $c_1$. Define
$C^\star_1$ as the ellipse obtained by applying the half-turn with
center $r$ to $C_1$. The ellipses $C_0$, $C_1$, and $C^\star_1$ are
described by matrices $M_0$, $M_1$, and $M^\star_1$ with respective
eigenvalues
\begin{equation*}
  e(M_0) = (1,\nu_{0,1},\nu_{0,2}),
  \quad
  e(M_1) = e(M^\star_1) = (1,\nu_{1,1},\nu_{1,2}).
\end{equation*}
We would like to make some admissible assumptions on these
eigenvalues. Because of $\area(C_0) = \area(C_1)$, we have
\begin{equation}
  \label{eq:10}
  1 < \nu_{0,1} \le \nu_{1,1} \le \nu_{1,2} \le \nu_{0,2}.
\end{equation}
If $\nu_{0,1} = \nu_{1,1}$ or $\nu_{1,2} = \nu_{0,2}$, \eqref{eq:10}
holds with equality throughout and both ellipses are actually
congruent circles. In this case a simple construction produces a
smaller enclosing circle (Figure~\ref{fig:two-circles}): Denote the
two intersection points of $C_0$ and $C_1$ by $s_0$ and $s_1$.  By
elementary hyperbolic geometry, the circle $S$ over the diameter
$s_0$, $s_1$ is smaller than $C_0$ and $C_1$ and it contains the
common interior of $C_0$ and~$C_1$.

Thus, the case of two congruent circles can be excluded and we may
assume that the eigenvalues of $M_0$ and $M_1$ are ordered
according to
\begin{equation}
  \label{eq:11}
  1 < \nu_{0,1} < \nu_{1,1} \le \nu_{1,2} < \nu_{0,2}.
\end{equation}

Now we compute the derivative of the area function \eqref{eq:6} with
respect to $\lambda$. For that purpose, we assume that $C_0$ is given
by the normal form~\eqref{eq:4} and $C_1$ is obtained from an ellipse
in this normal form by a hyperbolic rotation, that is,
\begin{equation}
  \label{eq:12}
  M_0 = \diag(-1, \nu_{0,1}, \nu_{0,2}),
  \quad
  M_1 = (Q^{-1})^\tp \cdot \diag(-1, \nu_{1,1}, \nu_{1,2}) \cdot Q^{-1}
\end{equation}
with the hyperbolic rotation matrix
\begin{equation}
  \label{eq:13}
  \begin{gathered}
    Q =
    \begin{pmatrix}
      q_0^2 + q_1^2 + q_2^2 + q_3^2 & 2(q_0 q_3 + q_1 q_2)          & 2(q_1 q_3 - q_0 q_2) \\
      2(q_0 q_3 - q_1 q_2)          & q_0^2 - q_1^2 - q_2^2 + q_3^2 & 2(q_0 q_1 - q_2 q_3) \\
      2(-q_0 q_2 - q_1 q_3)         & 2(-q_0 q_1 - q_2 q_3)         & q_0^2 - q_1^2 + q_2^2 - q_3^2
    \end{pmatrix},\\
    q_0^2 + q_1^2 - q_2^2 - q_3^2 = 1.
  \end{gathered}
\end{equation}
Up to the irrelevant sign of $q_1$ and an index shift, this is
precisely Equation~(5) of \cite{oezdemir06:_minkowski_rotations}.  The
rotation angle $-2\xi$ is given by $q_0 = \cos\xi$, the axis direction
is $(-q_1, q_2, q_3)^\tp$. A hyperbolic half-turn is obtained by
substituting $q_0 = 0$ into \eqref{eq:13}. In this case, $Q \cdot Q$
indeed equals the unit matrix $\diag(1,1,1)$.

The matrix $M_\lambda$ to $C_\lambda$ is computed according to
Definition~\ref{def:in-between-ellipse}. Its ordered eigenvalues
$(\nu_0, \nu_1, \nu_2)$ are functions of $\lambda$. In the vicinity of
$\lambda = 0$ we have $\nu_0(\lambda) > 0$ and $1 < \nu_1(\lambda) <
\nu_2(\lambda)$. These eigenvalues are implicitly defined as roots of
the characteristic polynomial $P(\lambda, \nu(\lambda)) =
\det(M_\lambda - \nu I)$ of $M_\lambda$ where $I$ is the matrix
defined in \eqref{eq:3}. For $\lambda = 0$ we know the values of
these roots:
\begin{equation*}
  \nu_0(0) = 1,\quad
  \nu_1(0) = \nu_{0,1},\quad
  \nu_2(0) = \nu_{0,2}.
\end{equation*}
By implicit derivation we have
\begin{equation*}
  \dod{\nu_i}{\lambda}(0) =
  -\frac{\pd{P}{\lambda}(0, \nu_i(0))}{\pd{P}{\nu}(0, \nu_i(0))},
  \quad i = 0,1,2.
\end{equation*}
Furthermore, we can compute the partial derivatives
\begin{equation*}
  \dpd{\area(\nu_0, \nu_1, \nu_2)}{\nu_i}, \quad i = 0,1,2
\end{equation*}
of \eqref{eq:6}. Using the chain rule
\begin{equation*}
  \dpd{\area(C_\lambda)}{\lambda}\Big|_{\lambda=0}
  = \dpd{\area}{\nu_0} \dpd{\nu_0}{\lambda}\Big|_{\lambda=0}
  + \dpd{\area}{\nu_1} \dpd{\nu_1}{\lambda}\Big|_{\lambda=0}
  + \dpd{\area}{\nu_2} \dpd{\nu_2}{\lambda}\Big|_{\lambda=0},
\end{equation*}
we find
\begin{equation}
  \label{eq:14}
  \dpd{\area(C_\lambda)}{\lambda}\Big|_{\lambda=0} =
  -\frac{1}{2}\int_{-\pi}^\pi \frac{D}{N} \dif\varphi
\end{equation}
where
\begin{equation*}
  \begin{aligned}
    D =\mbox{} &((q_{1,2}^2 \nu_{1,2} + q_{1,3}^2 \nu_{1,2} - q_{1,1}^2)\nu_{0,1} +
        q_{2,2}^2 \nu_{1,1} + q_{2,3}^2 \nu_{1,2} - q_{2,1}^2) \sin^2\varphi + \\
     &((q_{1,2}^2 \nu_{1,1} + q_{1,3}^2 \nu_{1,2} - q_{1,1}^2)\nu_{0,2} +
       q_{3,2}^2 \nu_{1,1} + q_{3,3}^2 \nu_{1,2} - q_{3,1}^2) \cos^2\varphi,\\
    N =\mbox{}& (\nu_{0,1} \sin^2 \varphi + \nu_{0,2} \cos^2 \varphi - 1)^{3/2}
  (\nu_{0,1} \sin^2 \varphi + \nu_{0,2} \cos^2 \varphi)^{1/2}
  \end{aligned}
\end{equation*}
and $q_{i,j}$ are the entries of the matrix~\eqref{eq:13}. It will be
convenient to write \eqref{eq:14} in terms of the first and second
complete elliptic integrals
\begin{equation}
  \label{eq:15}
  K(z) = \int_0^1 \frac{1}{\sqrt{1-t^2} \sqrt{1-z^2t^2}} \dif t
  \quad\text{and}\quad
  E(z) = \int_0^1 \frac{\sqrt{1-z^2t^2}}{\sqrt{1-t^2}} \dif t.
\end{equation}
Since we will evaluate them only at
\begin{equation*}
  f = \sqrt{\frac{\nu_{0,2}-\nu_{0,1}}{(\nu_{0,2} - 1) \nu_{0,1}}},
\end{equation*}
we use the abbreviations $\ellE := E(f)$ and $\ellK := K(f)$. By
\eqref{eq:11}, $f$ is always real and between 0 and 1. Substituting
\begin{equation*}
  \varphi = \arcsin\sqrt{\frac{\nu_{0,2}-x}{\nu_{0,2}-\nu_{0,1}}},
\end{equation*}
and noting that $2(\nu_{0,2}-x)(x-\nu_{0,1}) = 2(\nu_{0,1} -
\nu_{0,2}) \cos\varphi \sin\varphi$ we can express the derivative of
the area function in terms of $\ellE$ and $\ellK$:
\begin{multline}
  \label{eq:16}
  \dpd{\area(C_\lambda)}{\lambda} \Big|_{\lambda=0} =
  \frac{2}
       {\sqrt{(\nu_{0,2}-1)\nu_{0,1}}(\nu_{0,2}-\nu_{0,1})(\nu_{0,1}-1)}\\
  (
      \cA (q_{1,2}^2 \nu_{1,1} + q_{1,3}^2 \nu_{1,2} - q_{1,1}^2) +
      \cB (q_{2,2}^2 \nu_{1,1} + q_{2,3}^2 \nu_{1,2} - q_{2,1}^2) +
      \cC (q_{3,2}^2 \nu_{1,1} + q_{3,3}^2 \nu_{1,2} - q_{3,1}^2)
  )
\end{multline}
where
\begin{equation}
  \label{eq:17}
  \begin{gathered}
    \cA = -\nu_{0,1} (\nu_{0,2}-\nu_{0,1}) \ellE,
    \quad
    \cB = \nu_{0,2} (\nu_{0,1}-1) \ellK - \nu_{0,1} (\nu_{0,2}-1) \ellE,\\
    \cC = \nu_{0,1} (\nu_{0,1}-1) (\ellE-\ellK).
  \end{gathered}
\end{equation}

Having computed \eqref{eq:16}, the preparatory work for the final
(big) step in the proof of Theorem~\ref{th:1} is completed. We
formulate the last step as a lemma:

\begin{lemma}[Hyperbolic Half-Turn Lemma]
  \label{lem:half-turn}
  Consider three ellipses $C_0$, $C_1$, $C^\star_1$ of equal
  area. Assume that
  \begin{itemize}
  \item $C_0$ and $C^\star_1$ are concentric,
  \item $C_1$ is obtained from $C^\star_1$ by a half-turn,
  \item the eigenvalues $\nu_{i,1}$, $\nu_{i,2}$ of the normalized
    matrix $M_i$ to $C_i$ ($i=0,1$) satisfy ~\eqref{eq:11}, and
  \item $H(\nu_{i,1}, \nu_{i,2}) \le 0$ where $H$ is defined in
    Equation~\eqref{eq:1}.
  \end{itemize}
  Then the area of $C_\lambda = (1-\lambda)C_0 + \lambda C_1$ is
  smaller than the area of $C^\star_\lambda = (1-\lambda)C_0 + \lambda
  C^\star$, at least in the vicinity of $\lambda = 0$.
\end{lemma}

In order to proof Lemma~\ref{lem:half-turn}, we compare the
derivatives of the areas of $C_\lambda$ and $C^\star_\lambda$ with
respect to $\lambda$ at $\lambda = 0$. The ellipse $C^\star_1$ can be
obtained from an ellipse in normal form~\eqref{eq:4} by a rotation
about $(1, 0, 0)^\tp$ through $\zeta$. We can compute the matrix $M_1$
as in~\eqref{eq:12} by substituting
\begin{equation*}
  q_0 = \cos\tfrac{\zeta}{2},
  \quad
  q_1 = -\sin\tfrac{\zeta}{2},
  \quad
  q_2 = q_3 = 0
\end{equation*}
into the matrix~\eqref{eq:13}. Plugging this into
Equation~\eqref{eq:16} yields
\begin{equation*}
  \frac{1}{2}\dpd{\area(C^\star_\lambda)}{\lambda} \Big|_{\lambda=0} =
  \frac{D^\star_1}{N^\star_1}
\end{equation*}
where
\begin{equation*}
  \begin{aligned}
    D^\star_1 &= -\cA +
                 (\cB\cos^2\zeta+\cC\sin^2\zeta)\nu_{1,1} +
                 (\cB\sin^2\zeta+\cC\cos^2\zeta)\nu_{1,2},\\
    N^\star_1 &= \sqrt{(\nu_{0,2}-1) \nu_{0,1}}(\nu_{0,2} - \nu_{0,1})(\nu_{0,1} - 1),
  \end{aligned}
\end{equation*}
and $\cA$, $\cB$, $\cC$ are as in~\eqref{eq:17}.

The ellipse $C_1$ is obtained by a half-turn from $C^\star_1$ about
the rotation axis defined by the unit vector
$r=(r_1,r_2,r_3)^\tp$. The matrix $Q$ in~\eqref{eq:13} is the product
of the rotation matrix about $(1, 0, 0)^\tp$ through $\zeta$ and a
half-turn rotation matrix about the unit vector $r$. The later is
obtained by substituting
\begin{equation*}
  q_0 = 0,\quad
  q_1 = -r_1,\quad
  q_2 =  r_2,\quad
  q_3 =  r_3
\end{equation*}
into Equation~\eqref{eq:13}. Plugging the entries of the product
matrix into \eqref{eq:16} yields
\begin{equation*}
  \frac{1}{2} \dpd{\area(C_1)}{\lambda} \Big|_{\lambda=0} \equiv \frac{D_1}{N_1}
  \mod (r_1^2-r_2^2-r_3^2-1)
\end{equation*}
where $N_1 = N_1^\star$ and
\begin{equation*}
  \begin{aligned}
    D_1 =&\phantom{\mbox{}+\mbox{}}4r_2r_3((2\cA+\cB+\cC)r_1^2+(\cB-\cC)r_2^2-(\cB-\cC)r_3^2)(\nu_{1,1}-\nu_{1,2})\sin\zeta\cos\zeta\\
         &+ (4(\cA+\cB)r_1^2r_2^2-4(\cA+\cC)r_1^2r_3^2-8(\cB-\cC)r_2^2 r_3^2+\cB-\cC)(\nu_{1,1}-\nu_{1,2})\cos^2\zeta\\
         &+ 4(\cA+\cB)r_1^2r_2^2(\nu_{1,2}-1)+4(\cA+\cC)r_1^2r_3^2(\nu_{1,1}-1)\\
         &+ 4(\cB-\cC)r_2^2r_3^2(\nu_{1,1}-\nu_{1,2})+\cC\nu_{1,1}+\cB\nu_{1,2}-\cA.
  \end{aligned}
\end{equation*}
Now we are going to prove the inequality $D_1 - D^\star_1 < 0$ for
$\zeta \in [0,\frac{\pi}{2}]$. We substitute $\zeta = 2\arctan t$ into
its left-hand side and obtain a rational expression in $t$. Clearing
the positive denominator $(1+t^2)^2$, we are left with a polynomial
$P(t)$ of degree four whose negativity on $[0,1]$ has to be shown. To
do this, we write $P(t)$ with respect to the Bernstein basis as
\begin{equation}
  \label{eq:18}
  P(t) = \sum_{i=0}^4 p_i B^4_i(t)
  \quad\text{where}\quad
  B^4_i(t) = \binom{4}{i}(1-t)^{4-i}t^i
\end{equation}
and show non-positivity of the coefficients $p_1$, $p_2$, $p_3$ and
negativity of the coefficients $p_0$ and $p_4$.  After a
straightforward basis transformation and reducing modulo
$r_1^2-r_2^2-r_3^2-1$ we find
\begin{equation*}
  p_0 = 4(
     (\nu_{1,1}-1)(\cA+\cB)r_1^2r_2^2 +
     (\nu_{1,2}-1)(\cA+\cC)r_1^2r_3^2 +
     (\nu_{1,2}-\nu_{1,1})(\cB-\cC)r_2^2r_3^2
  ),
\end{equation*}
\begin{multline*}
  p_1 =
  4 (
     (\nu_{1,1}-1)(\cA+\cB)r_1^2r_2^2
     + (\nu_{1,2}-1)(\cA+\cC)r_1^2r_3^2
     + (\nu_{1,2}-\nu_{1,1})(\cB-\cC)r_2^2r_3^2
  )\\
  +2r_2r_3(\nu_{1,2}-\nu_{1,1})(
      -(2\cA+\cB+\cC)r_1^2
      -(\cB-\cC)r_2^2
      +(\cB-\cC)r_3^2
  ),
\end{multline*}
\begin{multline*}
  3p_2 =
  8(\nu_{1,1}+\nu_{1,2}-2)(
      (\cA+\cB)r_1^2r_2^2
     +(\cA+\cC)r_1^2r_3^2
  )\\
  +12r_2r_3(\nu_{1,2}-\nu_{1,1})(
     -(2\cA+\cB+\cC)r_1^2
     -(\cB-\cC)r_2^2
     +(\cB-\cC)r_3^2
  ),
\end{multline*}
\begin{multline*}
  p_3 =
  8(
      (\nu_{1,2}-1)(\cA+\cB)r_1^2r_2^2
      +(\nu_{1,1}-1)(\cA+\cC)r_1^2r_3^2
      -(\nu_{1,2}-\nu_{1,1})(\cB-\cC)r_2^2r_3^2
  )\\
  +4r_2r_3(\nu_{1,2}-\nu_{1,1})(
      -(2\cA+\cB+\cC)r_1^2
      -(\cB-\cC)r_2^2
      +(\cB-\cC)r_3^2
  ),
\end{multline*}
\begin{equation*}
  p_4 =
  16(
     (\nu_{1,2}-1) (\cA+\cB)r_1^2r_2^2
     + (\nu_{1,1}-1)(\cA+\cC)r_1^2r_3^2
     - (\nu_{1,2}-\nu_{1,1})(\cB-\cC)r_2^2r_3^2
  ).
\end{equation*}

Recall now Equation~\eqref{eq:11} ($1 < \nu_{0,1} < \nu_{1,1} \le
\nu_{1,2} < \nu_{0,2}$) and $r_1^2 - r_2^2 - r_3^2 = 1$ and observe
that
\begin{itemize}
\item $\cA < \cB < \cC$; this is proved in Lemma~\ref{lem:2} and
  Lemma~\ref{lem:3} in the appendix.
\item $\cC < 0$; this follows from $\ellE < \ellK$ and $\nu_{0,1} >
  1$.
\end{itemize}
Under these conditions, the negativity of $p_0$ is clear except when
$r_2 = r_3 = 0$. But this is the concentric case $C_1 = C^\star_1$ and
need not be considered. The non-positivity of the coefficients $p_1$,
$p_2$, $p_3$, and the negativity of $p_4$ is shown in
Lemmas~\ref{lem:6} and \ref{lem:5} below.  This concludes the proof of
the Half-Turn Lemma and, thus, also the proof of Theorem~\ref{th:1}.

\begin{example}
  \label{ex:1}
  We use the prerequisites of Theorem~\ref{th:1} on the radii $r$ and
  $R$ in Lemma~\ref{lem:6}.  But one might wonder whether the
  Half-Turn Lemma remains true without these assumptions.  The answer
  to this question is negative. We can provide and example, where the
  polynomial $P(t)$ attains positive values on $(0,1)$.

  Substituting $r_1^2 = r_2^2 + r_3^2 + 1$, the coefficient $p_1$ can
  be written as
  \begin{equation*}
    \begin{aligned}
      p_1 &= \nu_{1,2}(
           4r_2^2r_3(\cA+\cB)(r_3-r_2)
          +4r_3^2(\cA+\cC)(r_3^2-r_2r_3+1)
          -2(2\cA+\cB+\cC)r_2r_3
      ) \\
      &+\nu_{1,1}(
           4(\cA+\cB)r_2^2(r_2^2+r_2r_3+1)
          +4(\cA+\cC)r_2r_3^2(r_2+r_3)
          +2(2\cA+\cB+\cC)r_2r_3
      ) \\
      &-4(1+r_2^2+r_3^2)(r_2^2(\cA+\cB)+r_3^2(\cA+\cC)).
    \end{aligned}
  \end{equation*}
  Assuming $r_2$, $r_3 > 0$ we see that
  \begin{itemize}
  \item the coefficient of $\nu_{1,1}$ is always negative and
  \item it is possible to choose $\nu_{0,1}$, $\nu_{0,2}$, $r_2$, and
    $r_3$ so that the coefficient of $\nu_{1,2}$ is positive.
  \end{itemize}
  Consequently $p_1$ can be made positive for large $\nu_{1,2}$. The
  choice
  \begin{equation*}
    \nu_{0,1} = \nu_{1,1} = 1.1,\quad
    \nu_{0,2} = \nu_{1,2} = 90,\quad
    r_2 = 0.9 \cos(0.3),\quad
    r_3 = 0.9 \sin(0.3),
  \end{equation*}
  accomplishes this and even makes $P(t)$ attain positive values for
  $t \in (0,1)$ (the zeros are $t \approx 0.1272$ and $t \approx
  0.1389$).  Note that this does not imply
  \begin{equation*}
    \pd{\area(C_\lambda)}{\lambda}\Big|_{\lambda = 0} > 0
  \end{equation*}
  and, thus, constitutes no counter-example to the statement that the
  area of $C_\lambda$ is smaller than the area of $C_0$ and $C_1$. We
  are not aware of such a counter-example.
\end{example}

\section{Conclusion and future research}
\label{sec:conclusion}

We proved uniqueness results for minimal enclosing ellipses in the
hyperbolic plane. The general result (Theorem~\ref{th:1}) involves
rather cumbersome but straightforward calculations.  The differences
to the elliptic case are mainly in the final estimates for the
coefficients of the polynomial $P(t)$ in \eqref{eq:18} and can be
found in the appendix.

It is apparent that Theorem~\ref{th:1} leaves room for improvements.
Pushing back the frontier dictated by the inequality \eqref{eq:1} in
Theorem~1 would be nice. Substantial steps towards answering the
question whether the minimal area ellipse to all compact and
full-dimensional sets $F$ in the hyperbolic plane is unique or not
would be great.

Note that there is a subtle difference to the situation in the
elliptic plane. In \cite{weber10:_minimal_area_conics}, we presented
an example from which we inferred that uniqueness in the elliptic
plane cannot be proved by means of our construction of in-between
conics. In the hyperbolic plane, we are not aware of such a
configuration. Example~\ref{ex:1} only shows that the estimate of the
derivative of the area function is insufficient.  Thus, there is a
certain hope that a general uniqueness result can be proved by means
of our construction.

Since uniqueness or non-uniqueness of minimal enclosing ellipses in
the elliptic and hyperbolic plane remains a difficult topic, one might
try to aim at a weaker result and consider only ``typical'' (in the
sense of Baire categories, see \cite{gruber85:_results_baire,%
  gruber93:_baire_categories}) convex sets~$\conv{F}$.

We would also like to mention that \cite{weber10:_minimal_area_conics}
and this article are the only results on extremal quadrics in
non-Euclidean geometries that we are aware of.  We can conceive
numerous possibilities for generalizations.  They pertain to the
dimension of the surrounding space, the type of the quadric and
enclosed set (for example minimal enclosing hyperbolas to line sets as
in \cite{schroecker07:_minim_hyper}) the measure for the quadric's
size (volume, surface area etc.), and the replacement of ``minimal
enclosing'' by ``maximal inscribed'' quadrics.

The attentive reader will have noticed that our method of proving
Theorem~1 can be adapted to these generalizations.  Having defined an
``in-between'' quadric $Q_\lambda$ by means of a suitable matrix
convex combination, it might be infeasible to compute the size
$Q_\lambda$ in a form that allows further processing.  But, provided
the size function's derivatives with respect to the matrix eigenvalues
can be computed, it is, at least in principle, possible to obtain an
explicit formula for the derivative of the size of $Q_\lambda$ for
$\lambda = 0$.  Its negativity has to be shown so that the uniqueness
problem is made accessible to numerous tools and techniques related to
inequalities.

In the Euclidean setting, the mere uniqueness result is less important
than John's characterization of it via his famous decomposition of the
identity. The original reference is the old paper
\cite{john48:_studies_and_essays}. But, following
\cite{ball92:_ellipsoids_max_vol}, many contemporary authors
considered this topic \cite{gruber05:_arithmetic_proof,%
  gruber08:_john_type,%
  lutwak05:_john_ellipsoids,%
  bastero02:_johns_decomposition,%
  gordon04:_johns_decomposition}. Elliptic and hyperbolic versions of
John's characterization seem to be a worthy topic of future research.

% \appendix
\section*{Appendix. Proofs of auxiliary results}

\begin{lemma}
  \label{lem:2}
  For $\cA$ and $\cB$ as in~\eqref{eq:17} we have $\cA < \cB$.
\end{lemma}

\begin{proof}
  We show that $A-B < 0$. By \eqref{eq:17} we have
  \begin{equation*}
    \cA - \cB = (\nu_{0,1}-1) (\nu_{0,1} \ellE - \nu_{0,2} \ellK).
  \end{equation*}
  This is negative because of $\ellE < \ellK$ and $1 < \nu_{0,1} <
  \nu_{0,2}$.
\end{proof}

\begin{lemma}
  \label{lem:3}
  For $\cB$ and $\cC$ as in~\eqref{eq:17} we have $\cB < \cC$.
\end{lemma}

\begin{proof}
  We let $\cD = \cB - \cC$ and view $\cD$ as a function of $\nu_{0,1}$
  and $\nu_{0,2}$. Its negativity for $1 < \nu_{0,1} < \nu_{0,2}$
  follows from three facts:
  \begin{itemize}
  \item $\cD = 0$ for $\nu_{0,1} = \nu_{0,2}$ (this is obvious because
    in this case we have $\ellE = \ellK$),
  \item $\pd{\cD}{\nu_{0,2}} = 0$ for $\nu_{0,1} = \nu_{0,2}$, and
  \item $\cD$ is concave in $\nu_{0,2}$ for $1 < \nu_{0,1} <
    \nu_{0,2}$.
  \end{itemize}
  We compute the first partial derivative of $\cD$ with respect to
  $\nu_{0,2}$:
  \begin{equation*}
    \pd{\cD}{\nu_{0,2}} =
    \frac{(2\nu_{0,2}(1-\nu_{0,2})+\nu_{0,1}-1)\nu_{0,1}\ellE}
         {2\nu_{0,2}(\nu_{0,2}-1)} +
    \frac{(2\nu_{0,2}-1)(\nu_{0,1}-1)\ellK}
         {2(\nu_{0,2}-1)}.
  \end{equation*}
  It vanishes for $\nu_{0,1} = \nu_{0,2}$. The second partial
  derivative of $\cD$ with respect to $\nu_{0,2}$ equals
  \begin{equation*}
    \dpd[2]{\cD}{\nu_{0,2}} =
    \frac{(\nu_{0,1}-1)}{4\nu_{0,2}^2(\nu_{0,2}-1)^2} J_1
    \quad\text{where}\quad
      J_1 = \nu_{0,2}(\nu_{0,1}-1)\ellK-\nu_{0,1}(5\nu_{0,2}-2)\ellE.
  \end{equation*}
  We have to show that it is negative.  The factor before $J_1$ is
  positive. To see the negativity of $J_1$ itself we write it in the
  integral form (see~\eqref{eq:15})
  \begin{equation*}
    J_1 = \int_0^1 \frac{J_2}{\sqrt{1-t^2} \sqrt{1-f^2 t^2}} \dif t,
  \end{equation*}
  where
  \begin{equation*}
    J_2 = 
    \nu_{0,2}(\nu_{0,1}-1) - \nu_{0,1} (5\nu_{0,2}-2)
      \Bigl(
      1 - t^2 \frac{\nu_{0,2}-\nu_{0,1}}{\nu_{0,1} (\nu_{0,2}-1)}
      \Bigr).
  \end{equation*}
  The term $J_2$ is linear in $t^2$. For $t=0$ it equals $2\nu_{0,1}
  (1-2\nu_{0,2}) - \nu_{0,2} < 0$ and for $t=1$ it equals $-\nu_{0,2}
  (\nu_{0,1}-1) (4\nu_{0,2}-1)/(\nu_{0,2}-1) < 0$. Thus, $J_2 < 0$ for
  $t \in [0,1]$. This implies $J_1 < 0$ and we see that $\cD$ is
  indeed concave for $1 < \nu_{0,1} < \nu_{0,2}$.
\end{proof}

We will deduce non-positivity of the Bernstein coefficients $p_1$,
\ldots, $p_3$ from the inequality \eqref{eq:1} and the additional
inequalities
\begin{align}
  h_1(\nu_1,\nu_2) &:= \nu_2 - 5\nu_1 + 4 \le 0,\label{eq:19}\\
  h_2(\nu_1,\nu_2) &:= -5\nu_1^2+\nu_1\nu_2+\nu_1+\nu_2+2 \le 0,\label{eq:20}\\
  h_3(\nu_1,\nu_2) &:= \nu_2^2 -5\nu_1\nu_2-2\nu_1+4\nu_2+2 \le 0,\label{eq:21}\\
  h_4(\nu_1,\nu_2) &:= 5\nu_2^2 - 13\nu_1\nu_2 - 2\nu_1 + 6\nu_2 + 4 \le 0,\label{eq:22}\\
  h_5(\nu_1,\nu_2) &:= -5\nu_1^2+\nu_1\nu_2-\nu_1+3\nu_2+2 \le 0,\label{eq:23}\\
  h_6(\nu_1,\nu_2) &:= \nu_2^2-5\nu_1\nu_2+2\nu_2+2 \le 0,\label{eq:24}
\end{align}
which are all simple consequences of \eqref{eq:1}.  We state this in
Lemma~\ref{lem:9}, below.

The assumptions of Theorem~\ref{th:1} guarantee that these
inequalities are fulfilled for $\nu_1 = \nu_{0,1}$, $\nu_2 =
\nu_{0,2}$. Thus, we only have to show that the inequalities are
satisfied on the set
\begin{equation*}
  U := \{(\nu_1,\nu_2) \mid 1 < \nu_1 < \nu_2\}.
\end{equation*}

\begin{lemma}
  \label{lem:9}
  {\normalfont (a)} If a point $(\nu_1,\nu_2) \in U$ satisfies the
  inequality \eqref{eq:1}, it also satisfies the inequalities
  \eqref{eq:19}--\eqref{eq:24}.\par
  {\normalfont (b)} If a point $(\nu_1^\star, \nu_2^\star)$ satisfies
  the inequalities \eqref{eq:1}, and \eqref{eq:19}--\eqref{eq:24}, the
  same is true for all points $(\nu_1,\nu_2) \in U$ with $\nu_1 \ge
  \nu_1^\star$ and $\nu_2 \le \nu_2^\star$.
\end{lemma}

\begin{proof}
  It is an elementary exercise to verify that the curves defined by
  the implicit equations $H$ and $h_j$ over the closure of $U$ are the
  graphs of strictly monotone increasing functions of $F(\nu_1)$ and
  $f_j(\nu_1)$
  \begin{equation*}
    F, f_j\colon [1,\infty) \to [1,\infty),
  \end{equation*}
  for $j \in \{1,\ldots,6\}$. This implies assertion (b).

  The functions $f_j-F$ are strictly monotone increasing as well.
  This and the observation $F(1) = f_j(1) = 1$ implies assertion~(a).
\end{proof}

Figure~\ref{fig:hyp-graph} displays the hyperbolas $H(\nu_1,\nu_2) =
0$ and, as an example, $h_2(\nu_1,\nu_2) = 0$ together with the line
$\nu_1 = \nu_2$. The remaining curves are depicted in light-gray. The
region $U$ is dotted.

\begin{lemma}
  \label{lem:6}
  The coefficients $p_1$, $p_2$, and $p_3$ are not positive.
\end{lemma}

\begin{proof}
  % $p_1$ is of the form $p_0 + 2 r_2 r_3 (\nu_{1,2} - \nu_{1,1})
  % p^*$. In Lemma~\ref{lem:4} we show that $p^*$ is positive. So
  % in the case of $r_2 r_3 \leq 0$ the proof is finished. Therefore we
  % assume $r_2 r_3 > 0$.
  
  We substitute $r_1^2 = 1 + r_2^2 + r_3^2$ into $p_1$ and observe
  that $p_1 = 0$ and $\pd{p_1}{r_2} = 0$ if $r_2 = r_3 = 0$.  The
  lemma's claim holds true if we can show that the $r_2$-parameter
  lines of $p_1$, viewed as a function of $r_2$ and $r_3$, are
  strictly concave, that is,
  \begin{multline}
    \label{eq:25}
    \dpd[2]{p_1}{r_2} = 8(\cA+\cB)(\nu_{1,1}-1)(6 r_2^2+1)\\
    +8((\cA+\cB)\nu_{1,2} +(\cA+\cC)\nu_{1,1} -(2\cA+\cB+\cC))r_3^2
    -24(\nu_{1,2}-\nu_{1,1})(\cA+\cB) r_2 r_3 < 0.
  \end{multline}
  The coefficient of $r_2r_3$ is positive, the remaining terms are
  negative.  By the inequality of arithmetic and geometric means we
  have $r_2r_3 \le (r_2^2+r_3^2)/2$. We insert this into \eqref{eq:25}
  to obtain
  \begin{multline}
    \label{eq:26}
    \dpd[2]{p_1}{r_2} \le
    4(\cA+\cB)(3(-\nu_{1,2}+5\nu_{1,1}-4)r_2^2 + 2(\nu_{1,1}-1))\\
    + 4(-(\cA+\cB)\nu_{1,2} +(5\cA+3\cB+2\cC)\nu_{1,1}-2(2\cA+\cB+\cC)) r_3^2 < 0.
  \end{multline}
  The first term is negative if $\nu_{1,2} -5\nu_{1,1} + 4 \le
  0$. This is implied by $\nu_{0,2} - 5\nu_{0,1} + 4 \le 0$ and thus
  follows from \eqref{eq:19}.  In the second term the coefficient of
  $r_3^2$ needs closer investigation. We want to show its
  negativity. By \eqref{eq:11} we have
  \begin{multline}
    \label{eq:27}
    (5\cA+3\cB+2\cC)\nu_{1,1} - (\cA+\cB)\nu_{1,2} - 2(2\cA+\cB+\cC) \leq \\
    (5\cA+3\cB+2\cC)\nu_{0,1} - (\cA+\cB)\nu_{0,2} - 2(2\cA+\cB+\cC).
  \end{multline}
  Using \eqref{eq:15} and \eqref{eq:17}, we write the term on the
  right in its integral form:
  \begin{multline}
    \label{eq:28}
    (5\cA+3\cB+2\cC)\nu_{0,1} - (\cA+\cB)\nu_{0,2} -2(2\cA+\cB+\cC)\\
    = (\nu_{0,2}-\nu_{0,1}) \int_0^1 \frac{J_3}{\sqrt{1-t^2}\sqrt{1-f^2
        t^2}}
  \end{multline}
  where
  \begin{equation*}
    J_3 = 
      -\frac{\nu_{0,2}-\nu_{0,1}}{\nu_{0,2}-1}(2\nu_{0,2}-7\nu_{0,1}+5)t^2
      +\nu_{0,2}(\nu_{0,1}+1)+\nu_{0,1}(1-5\nu_{0,1})+2.
  \end{equation*}
  We see that $J_3$ is linear in $t^2$. For $t=0$ and $t=1$ it attains
  the respective values
  \begin{align}
    J_3\big\vert_{t=0} & = -5\nu_{0,1}^2 + \nu_{0,1}\nu_{0,2} + \nu_{0,1} + \nu_{0,2} + 2,\label{eq:29} \\
    J_3\big\vert_{t=1} & = \frac{\nu_{0,1}-1}{\nu_{0,2}-1}(\nu_{0,2}^2-5\nu_{0,1}\nu_{0,2}-2\nu_{0,1}+4\nu_{0,2}+2).\label{eq:30}
  \end{align}
  The right-hand side of \eqref{eq:29} is not positive by
  \eqref{eq:20}. The right-hand side of \eqref{eq:30} is not positive
  by \eqref{eq:21}. We conclude that the integrand $J_3$ is not
  positive for $t \in [0,1]$ and the same is true for
  $\pd[2]{p_1}{r_2}$. Hence, the coefficient $p_1$ as a function of
  $r_2$ is concave with the maximum, $p_1 = 0$ attained at $r_2 = r_3
  = 0$. Thus, $p_1$ is not positive.

  The proofs of non-positivity of $p_2$ and $p_3$ run along exactly
  the same lines.  We only provide the relevant formulas and reduce
  the explanatory text between them to a minimum.
  Equations~\eqref{eq:25} and \eqref{eq:26} become
  \begin{equation*}
    \begin{aligned}
      3\dpd[2]{p_2}{r_2} &= 16(\cA+\cB)(\nu_{1,2}+\nu_{1,1}-2)(6r_2^2+1) +\\
      & \qquad 16(2\cA+\cB+\cC)(\nu_{1,2}+\nu_{1,1}-2)r_3^2 - 144(\cA+\cB)(\nu_{1,2}-\nu_{1,1})r_2r_3\\
      &\leq 8(\cA+\cB)(3(\nu_{1,2}+7\nu_{1,1}-8)r_2^2 +2(\nu_{1,2}+\nu_{1,1}-2)) +\\
      & \qquad 8(-(5\cA+7\cB-2\cC)\nu_{1,2} + (13\cA+11\cB+2\cC)\nu_{1,1} - 4(2\cA+\cB+\cC))r_3^2.
    \end{aligned}
  \end{equation*}
  Instead of \eqref{eq:27} and \eqref{eq:28} we have
  \begin{equation*}
    \begin{gathered}
      -(5\cA+7\cB-2\cC)\nu_{1,2} + (13\cA+11\cB+2\cC)\nu_{1,1} - 4(2\cA+\cB+\cC) \le \\
      -(5\cA+7\cB-2\cC)\nu_{0,2} + (13\cA+11\cB+2\cC)\nu_{0,1} - 4(2\cA+\cB+\cC) = \\
      (\nu_{0,2}-\nu_{0,1}) \int_0^1 \frac{J_4}{\sqrt{1-t^2}\sqrt{1-f^2 t^2}} \dif t
    \end{gathered}
  \end{equation*}
  where
  \begin{equation*}
    J_4 = -\frac{\nu_{0,2}-\nu_{0,1}}{\nu_{0,2}-1}(3(4\nu_{0,2}-5\nu_{0,1}+1)t^2
          + \nu_{0,2}(5\nu_{0,1}+7)+\nu_{0,1}(-13\nu_{0,1}-3)+4).
  \end{equation*}
  The non-positivity of $p_2$ follows from
  \begin{align*}
      J_4\big\vert_{t=0} & = -13\nu_{0,1}^2+5\nu_{0,1}\nu_{0,2}-3\nu_{0,1}+7\nu_{0,2}+4,\\
      J_4\big\vert_{t=1} & = \frac{\nu_{0,1}-1}{\nu_{0,2}-1}(5\nu_{0,2}^2-13\nu_{0,1}\nu_{0,2}-2\nu_{0,1}+6\nu_{0,2}+4),
  \end{align*}
  \eqref{eq:1} and \eqref{eq:22}.

  As to the coefficient $p_3$, Equations~\eqref{eq:25} and
  \eqref{eq:26} are replaced by
  \begin{equation*}
    \begin{aligned}
      \dpd[2]{p_3}{r_2} & = 16(\nu_{1,2}-1)(\cA+\cB)(6r_2^2+1) +                                     \\
                        & 16((\cA+\cC)\nu_{1,2} + (\cA+\cB)\nu_{1,1} - (2\cA+\cB+\cC))r_3^2 - 
                        48(\nu_{1,2}-\nu_{1,1})(\cA+\cB)r_2r_3                              \\
                        & \leq 8(\cA+\cB)(3(3\nu_{1,2}+\nu_{1,1}-4)r_2^2 +2(\nu_{1,2}-1)) +          \\
                        & \qquad 8(-(\cA+3\cB-2\cC)\nu_{1,2} +5(\cA+\cB)\nu_{1,1} -2(2\cA+\cB+\cC))r_3^2
    \end{aligned}
  \end{equation*}
  and \eqref{eq:27} and \eqref{eq:28} by
  \begin{equation*}
    \begin{gathered}
      -(\cA+3\cB-2\cC)\nu_{1,2} +5(\cA+\cB)\nu_{1,1} -2(2\cA+\cB+\cC) \le \\ 
      -(\cA+3\cB-2\cC)\nu_{0,2} +5(\cA+\cB)\nu_{0,1} -2(2\cA+\cB+\cC) = \\ 
      (\nu_{0,2}-\nu_{0,1}) \int_0^1 \frac{J_5}{\sqrt{1-t^2}\sqrt{1-f^2t^2}}
    \end{gathered}
  \end{equation*}
  where
  \begin{equation*}
    J_5 =
        -\frac{\nu_{0,2}-\nu_{0,1}}{\nu_{0,2}-1}((4\nu_{0,2}-5\nu_{0,1}+1)t^2
      +\nu_{0,2}(\nu_{0,1}+3) -\nu_{0,1}(5\nu_{0,1}+1) +2).
  \end{equation*}
  The non-positivity of $p_3$ follows from
  \begin{align*}
    J_5\big\vert_{t=0}& = -5\nu_{0,1}^2+\nu_{0,1}\nu_{0,2}-\nu_{0,1}+3\nu_{0,2}+2,\\
    J_5\big\vert_{t=1}& = \frac{\nu_{0,1}-1}{\nu_{0,2}-1}(\nu_{0,2}^2 -5\nu_{0,1}\nu_{0,2}+2\nu_{0,2}+2).
  \end{align*}
  and \eqref{eq:23},~\eqref{eq:24}.
\end{proof}

% \begin{lemma}
%   \label{lem:4}
%   The expression $-(2\cA+\cB+\cC)r_1^2 - (\cB-\cC)r_2^2 +
%   (\cB-\cC)r_3^2$ is positive.
% \end{lemma}

% \begin{proof}
%   Because of $\cA < \cB < \cC < 0$ we also have $2\cA+\cB+\cC < \cA <
%   \cB < \cB-\cC$. Moreover, $r_1^2-r_2^2-r_3^2 = 1$ implies $r_1^2 \ge
%   r_3^2+1 > r_3^2$. Thus, $(2\cA+\cB+\cC)r_1^2 - (\cB-\cC)r_3^2 < 0$
%   and, because of $(\cB-\cC)r_2^2 < 0$,
%   \begin{equation}
%     (2\cA+\cB+\cC)r_1^2 + (\cB-\cC)r_2^2 - (\cB-\cC)r_3^2 < 0.
%   \end{equation}
%   This is equivalent to this lemma's claim.
% \end{proof}

The negativity of the only remaining Bernstein coefficient can be
shown directly without resorting to Lemma~\ref{lem:9}:

\begin{lemma}
  \label{lem:5}
  The coefficient $p_4$ is negative.
\end{lemma}

\begin{proof}
  We can write
  \begin{equation}
    \label{eq:31}
    \frac{p_4}{16} = 
      (\nu_{1,2}-1)r_2^2((\cA+\cB) r_1^2 - (\cB-\cC)r_3^2)
    + (\nu_{1,1}-1)r_3^2((\cA+\cC) r_1^2 + (\cB-\cC)r_2^2).
  \end{equation}
  The proof is finished, if we can show that the coefficients of
  $(\nu_{1,2}-1) r_2^2$ and $(\nu_{1,1}-1) r_3^2$ in \eqref{eq:31} are
  negative. For the coefficient of $(\nu_{1,2}-1)r_2^2$ we argue as
  follows: $\cA < \cB < \cC < 0$ implies $\cA+\cB < \cA < \cB <
  \cB-\cC$ and $r_1^2 - r_2^2 - r_3^2 = 1$ implies $r_1^2 >
  r_3^2$. Thus, $(\cA+\cB)r_1^2-(\cB-\cC)r_3^2 < 0$.  The negativity
  of the coefficient of $(\nu_{1,2}-1)r_3^2$ follows from $\cA + \cC <
  0$ and $\cB-\cC < 0$.
\end{proof}

\section*{Acknowledgments}

The authors gratefully acknowledge support of this research by the
Austrian Science Foundation FWF under grant P21032 (Uniqueness Results
for Extremal Quadrics).

\bibliographystyle{plainnat}
\bibliography{mrabbrev,exquadb}

\end{document}